\documentclass{gtart}

\usepackage{rlepsf, amssymb, amsmath}

\input gtoutput
\volumenumber{3}\papernumber{11}\volumeyear{1999}
\pagenumbers{253}{268}\published{5 September 1999}
\proposed{Robion Kirby}
\seconded{Yasha Eliashberg, Tomasz Mrowka}
\received{16 June 1999}
\accepted{27 August 1999}

\newtheorem{thm}{Theorem}[section]
\newtheorem{lem}[thm]{Lemma}

\newtheorem{prop}[thm]{Proposition}
\newtheorem{rem}[thm]{Remark}
\newtheorem*{ttm}{Theorem}
\newtheorem{quest}{Question}

\def\C{\hbox{$\mathbb C$} }
\def\Z{\hbox{$\mathbb Z$} }

\def\R{\hbox{$\mathbb R$} }

\def\dfn#1{{\em #1}}

\def\co{\colon\thinspace}

\def\eqn#1{Equation~\eqref{#1}}
\def\tm#1{Theorem~\ref{#1}}
\def\lm#1{Lemma~\ref{#1}}

\begin{document}

%%%%%%%%%%%%%%%%%%%%%%%%%%%%%%%%%%%%%%%%%%%%%%%%%%%%%%%%%%
\title{Transversal torus knots}

\author{John B Etnyre}
\address{Mathematics Department, Stanford University\\Stanford, CA
94305, USA}
\email{etnyre@math.stanford.edu}
\url{http://math.stanford.edu/\char126 etnyre}

\begin{abstract}
We classify positive transversal torus knots in tight contact
structures up to transversal isotopy.
\end{abstract}

\asciiabstract{We classify positive transversal torus knots in tight contact structures up to transversal isotopy.}

\keywords{Tight, contact structure, transversal knots, torus knots}

\primaryclass{57M50, 57M25}
\secondaryclass{53C15}

\maketitlepage

% ************************************************************************
\section{Introduction}
% ************************************************************************

The study of special knots in contact three manifolds provided great insight
into the geometry and topology of three manifolds.  In particular, the study
of Legendrian knots (ones tangent to the contact planes) has been useful
in distinguishing homotopic contact structures on $T^3$  \cite{k} and
homology spheres \cite{am}.  Moreover, Rudolph \cite{r} has shown
that invariants of Legendrian knots can be useful in understanding
slicing properties of knots.  The first example of the use of knot 
theory in contact topology was in the work of Bennequin. In \cite{be}
Bennequin used transversal knots (ones transversal to the contact
planes) to show that $\R^3$ has exotic contact structures.  This was the
genesis of Eliashberg's insightful tight versus overtwisted dichotomy in three dimensional
contact geometry.

In addition to its importance in the understanding of contact geometry, the study
of transversal and Legendrian knots is quite interesting in its own right.
Questions concerning transversal and Legendrian knots have most prominently appeared
in \cite{el:knots} and Kirby's problem list \cite{kirby}.
Currently there are very few general theorems concerning the classification
of these knots. In \cite{el:knots}, Eliashberg classified transversal unknots
in terms of their self-linking number.  In \cite{ef}, Legendrian unknots
were similarly classified.  In this paper we will extend this classification 
to positive transversal torus knots\footnote{By ``positive transversal torus knot'' we mean a positive (right handed) 
	torus
	knot that is transversal to a contact structure.}.  In particular we prove:
\begin{ttm}
	Positive transversal 
	torus knots
	are transversely isotopic if and only if they have the same topological knot
	type and the same self-linking number.
\end{ttm}
In the process of proving this result we will examine {\em transversal stabilization}.
This is a simple method for creating one transversal knot from another.  By 
showing that all positive transversal torus knots whose self-linking number is less than 
maximal come from this stabilization process we are able to reduce the above theorem
to the classification of positive transversal torus knots with maximal self-linking number.
Stabilization also provides a general way to approach the classification problem for 
other knot types.  For example, we can reprove Eliashberg's classification
of transversal unknots using stabilization ideas and basic contact topology.

It is widely believed that the
self-linking number is not a complete invariant for transversal knots. However,
as of the writing of this paper, there is no known knot type whose transversal realizations
are not determined by their self-linking number. For Legendrian knots, in contrast,
Eliashberg and Hofer (currently unpublished)
and Chekanov \cite{chek} have produced examples of Legendrian knots that are not determined
by their corresponding invariants.

In Section~\ref{basic} we review some standard facts concerning contact geometry
on three manifolds.  In Section~\ref{sec:main} we prove our main theorem modulo some
details concerning the characteristic foliations on tori which are proved in
Section~\ref{tori} and some results on stabilizations proved in Section~\ref{stab}.
In the last section we discuss some open questions.

\rk{Acknowledgments}
The author gratefully acknowledges the support of an NSF 
Post-Doctoral Fellowship (DMS--9705949) and Stanford University.
Conversations with Y Eliashberg and E Giroux were helpful
in preparing this paper.

% ************************************************************************
\section{Contact structures in three dimensions}\label{basic}
% ************************************************************************

We begin by recalling some basic facts from contact topology.  For a more detailed
introduction, see \cite{a:etall, gi:survey}.
Recall an orientable plane field $\xi$ is a \dfn{contact structure} on a three manifold
if $\xi=\hbox{ker } \alpha$ where $\alpha$ is a nondegenerate 1--form for which
$\alpha\wedge d\alpha\not=0$.  Note $d\alpha$ induces an orientation on $\xi$.
Two contact structures are called \dfn{contactomorphic} if there is a diffeomorphism
taking one of the plane fields to the other.
A contact structure $\xi$ induces a singular
foliation on a surface $\Sigma$  by integrating the singular line field
$\xi\cap T\Sigma$. This is called the \dfn{characteristic foliation} and is denoted
$\Sigma_\xi$. Generically, the singularities are elliptic (if local degree is 1)
or hyperbolic (if the local degree is $-1$).  If $\Sigma$ is oriented then the singularities
also have a sign.  A singularity is positive (respectively negative) if the orientations
on $\xi$ and $T\Sigma$ agree (respectively disagree) at the singularity. 

\begin{lem}[Elimination Lemma \cite{gi:convex}]\label{lem:elimination}
	Let $\Sigma$ be a surface in a contact $3$--manifold $(M,\xi)$.  Assume
	that $p$ is an elliptic and $q$ is a hyperbolic singular point in 
	$\Sigma_\xi$, they both have the same sign and there is a leaf $\gamma$
	in the characteristic foliation $\Sigma_\xi$ that connects $p$ to $q$.
	Then there is a $C^0$--small isotopy $\phi\co \Sigma\times[0,1]\to M$ such that
	$\phi_0$ is the inclusion map, $\phi_t$ is fixed on $\gamma$ and outside any 
	(arbitrarily small) pre-assigned
	neighborhood $U$ of $\gamma$ and $\Sigma'=\phi_1(\Sigma)$ has no singularities
	inside $U$.
\end{lem}

It is important to note that after the above cancellation there is a curve in the 
characteristic foliation on which the singularities had previously sat. 
In the case of positive
singularities this curve will consist of the (closure of the) stable manifolds of the hyperbolic point
and {\em any} arc leaving the elliptic point (see \cite{ef, e:dis}), and similarly for the negative singularity case.
One may also reverse this process and add a canceling pair of singularities along 
a leaf in the characteristic foliation.  It is also important to note:

\begin{lem}
	The germ of the contact structure $\xi$ along a surface $\Sigma$ is determined 
	by $\Sigma_\xi$.
\end{lem}

Now recall that a contact structure $\xi$ on $M$ is called \dfn{tight} if no disk embedded in 
$M$ contains a limit cycle in its characteristic foliation, otherwise it is called 
\dfn{overtwisted}. The standard contact structure on $S^3$,
induced from the complex tangencies to $S^3=\partial B^4$ where $B^4$ is the
unit 4--ball in $\C^2$, is tight.

A closed curve $\gamma\co S^{1}\to M$ in a contact manifold $(M,\xi)$ is called 
\dfn{transversal} if $\gamma'(t)$ is transverse to $\xi_{\gamma(t)}$ 
for all $t\in S^{1}$.  Notice a transversal curve can be 
\dfn{positive} or \dfn{negative} according as $\gamma'(t)$ agrees with 
the co-orientation of $\xi$ or not. We will restrict our attention to 
positive transversal knots (thus in this paper ``transversal'' means 
``positive transversal''). It can be shown that any curve can be made 
transversal by a $C^{0}$ small isotopy. It will be useful to note:

\begin{lem}[See \cite{el:knots}]\label{extend}
	If $\psi_t\co S^1\to M$ is a transversal isotopy, then there
	is a contact isotopy $f_t\co M\to M$ such that $f_t\circ \psi_0=\psi_t$.
\end{lem}

Given a transverse knot $\gamma$ in $(M,\xi)$ that bounds a surface $\Sigma$
we define the \dfn{self-linking number}, $l(\gamma)$, of $\gamma$ as follows: 
take a nonvanishing vector field $v$ in $\xi\vert_{\gamma}$ that 
extends to a nonvanishing vector field in $\xi\vert_{\Sigma}$ and let $\gamma'$
be $\gamma$ slightly pushed along $v$. Define
$$l(\gamma,\Sigma)=I(\gamma',\Sigma),$$
where $I(\,\cdot\,,\,\cdot\,)$ is the oriented intersection number.
There is a nice relationship between $l(\gamma,\Sigma)$ and the 
singularities of the characteristic foliation of $\Sigma$.  Let 
$d_{\pm}=e_{\pm}-h_{\pm}$ where $e_{\pm}$ and $h_{\pm}$ are the 
number of $\pm$ elliptic and hyperbolic points in the characteristic 
foliation $\Sigma_{\xi}$ of $\Sigma$, respectively.  In \cite{be} it 
was shown that
\begin{equation}
	l=d_{-}-d_{+}.
\end{equation}
When $\xi$ is a {\em tight} contact structure and $\Sigma$ is a 
{\em disk},
Eliashberg \cite{el:twenty} has shown, using the elimination lemma,  how to 
eliminate all the positive hyperbolic and negative elliptic points 
from $\Sigma_{\xi}$. Thus in a 
tight contact structure when $\gamma$ is an unknot $l(\gamma,\Sigma)$ is always negative.
More generally one can show (see \cite{be, el:twenty}) that
\begin{equation}\label{lbound}
	l(\gamma)\leq -\chi(\Sigma),
\end{equation}
where $\Sigma$ is a Seifert surface for $\gamma$ and $\chi(\Sigma)$ is its Euler number.

Any odd negative integer can be realized as the self-linking number for some 
transversal unknot.
The first general result concerning the classification of transversal knots
was the following:

\begin{thm}[Eliashberg \cite{el:knots}]\label{unknots}
	Two transversal unknots are transversely isotopic if and only if 
	they have the same self-linking number.
\end{thm}

Let $\mathcal{T}$ be the transversal isotopy classes of transversal knots in $S^3$ with its unique
tight contact structure. Let $\mathcal{K}$ be the isotopy classes of knots in $S^3$.
Given a transversal knot $\gamma\in \mathcal{T}$ we have two pieces of information: its knot type 
$[\gamma]\in\mathcal{K}$ and its self-linking number $l(\gamma)\in\Z$.  Define 
\begin{equation}\label{tmap}
	\phi\co \mathcal{T}\to\mathcal{K}\times \Z : \gamma\mapsto ([\gamma],l(\gamma)).
\end{equation}
The main questions concerning transversal knots can be phrased in terms of the image
of this map and preimages of points.  In particular the above results say that $\phi$
is onto 
$$U=[\hbox{unknot}]\times\{\hbox{negative odd integers}\}$$ and $\phi$ is one-to-one
on $\phi^{-1}(U)$.

We will also need to consider Legendrian knots.  A knot $\gamma$ is a \dfn{Legendrian knot} if
it is tangent to $\xi$.  The contact structure $\xi$ defines a canonical framing on a Legendrian
knot $\gamma$. If $\gamma$ is null homologous we may associate a number to this framing which
we call the \dfn{Thurston--Bennequin invariant} of $\gamma$ and denote it $\hbox{tb}(\gamma)$.
If we let $\Sigma$ be the surface exhibiting the null homology of $\gamma$ then we may trivialize
$\xi$ over $\Sigma$ and use this trivialization to measure the rotation of $\gamma'(t)$ around
$\gamma$.  This number $r(\gamma)$ is called the \dfn{rotation number} of $\gamma$. Note that
the rotation number depends on an orientation on $\gamma$. From an oriented 
Legendrian knot $\gamma$ one can obtain canonical positive and negative transversal knots
$\gamma_\pm$ by pushing $\gamma$ by vector fields tangent to $\xi$ but transverse to
$\gamma'(t)$.  One may compute
\begin{equation}
	l(\gamma_\pm)=\hbox{tb}(\gamma)\mp r(\gamma).
\end{equation}
This observation combined with \eqn{lbound} implies
\begin{equation}\label{tb-bound}
	\hbox{tb}(\gamma)+|r(\gamma)|\leq -\chi(\Sigma).
\end{equation}
Consider an oriented (nonsingular) foliation $\mathcal{F}$ on a torus $T$. The foliation is 
said to have a \dfn{Reeb component} if two oppositely oriented periodic orbits cobound an
annulus containing no other periodic orbits. 

\begin{lem}\label{curveinT}
	Consider a torus $T$ in a contact three manifold $(M,\xi)$.  If the characteristic foliation
	on $T$ is nonsingular and contains no Reeb components then  
	any closed curve on $T$ may be isotoped to be transversal to $T_\xi$ or into a leaf
	of $T_\xi$. Moreover there is at most one homology class in $H_1(T)$ that can
	be realized by a leaf of $T_\xi$. 
\end{lem}

Now let $\xi$ be a tight contact structure on a solid torus $S$ with nonsingular characteristic
foliation on it boundary $T=\partial S$. It is easy to arrange for $T_\xi$ to have no Reeb 
components \cite{ml}.
Since $\xi$ is tight the lemma above implies the meridian $\mu$
can be made transversal to $T_\xi$. We say $S$ has self-linking number $l$ if $l=l(\mu)$
(ie, the self-linking number of $S$ is the self-linking number of its meridian).

\begin{thm}[Makar--Limanov \cite{ml}]\label{solid_tori}
	Any two tight contact structures on $S$ which induce the same nonsingular
	foliation on the boundary and have self-linking number $-1$ are contactomorphic.
\end{thm}

% ************************************************************************
\section{Positive transversal torus knots}\label{sec:main}
% ************************************************************************

Let $U$ be an unknot in a 3--manifold $M$, $D$ an embedded disk that it bounds
and $V$ a tubular neighborhood of $U$.  The boundary $T$ of $V$ is an embedded torus in $M$,
we call such a torus a \dfn{standardly embedded torus}. 
Let $\mu$ be the unique curve on $T$ that bounds a disk in $V$ and
$\lambda=D\cap V$.  Orient $\mu$ arbitrarily and 
then orient $\lambda$ so that $\mu, \lambda$ form a positive basis for $H_1(T)$ where 
$T$ is oriented as the boundary of $V$. Up to homotopy any curve in 
$T$ can be written as $p\mu + q\lambda$, we shall denote this curve by $K_{(p,q)}$.   
If $p$ and $q$ are relatively prime
then $K_{(p,q)}$ is called a {\em $(p,q)$--torus knot.}  If $pq>0$ we say $K(p,q)$ is a 
\dfn{positive} torus knot otherwise we call it \dfn{negative}. One
may easily compute that the Seifert surface of minimal genus for $K_{(p,q)}$
has Euler number $|p|+|q|-|pq|$.  Thus for a transversal torus knot Equation~\ref{lbound} implies 
\begin{equation}
	l(K_{(p,q)})\leq -|p|-|q|+|pq|.
\end{equation}
In fact, if $\overline{l}_{(p,q)}$ denotes the maximal self-linking number for a 
transversal $K_{(p,q)}$ then one may easily check that
\begin{equation}
	\overline{l}_{(p,q)}=-p-q+pq,
\end{equation}
if $p,q>0$, ie, for a positive torus knot. (Note: for a positive transversal torus knot 
Lemma~\ref{basis}
says we have $p,q>0$ not just $pq>0$.) From the symmetries involved in the definition of
a torus knot we may assume that $p>q$, which we do throughout the rest of the paper.
We now state our main theorem.

\begin{thm}\label{main}
	Positive transversal torus knots in a tight contact structure
	are determined up to transversal isotopy by 
	their knot type and their self-linking number.
\end{thm}

\begin{rem}
	{\em We may restate this theorem by saying
	the map $\phi$ defined in equation \ref{tmap} is one-to-one when restricted to
	$$(\hbox{pr}\circ \phi)^{-1}(\hbox{positive torus knots})$$ 
	(here $\hbox{pr}\co \mathcal{K}\times \Z\to \mathcal{K}$ is projection). 
	Moreover, the image of $\phi$ restricted to the above set is 
	$G=\cup_{(p,q)} K_{(p,q)}\times N(p,q)$ where the union
	is taken over relatively prime positive $p$ and $q$, and 
	$N(p,q)$ is the set of odd integers less than or equal to $-p-q+pq$.}
\end{rem}

We first prove the following auxiliary result:

\begin{prop}\label{aux}
	Two positive transversal $(p,q)$--torus knots $K$ and $K'$ in a tight contact
	structure with maximal self-linking number (ie, $l(K)=l(K')=\overline{l}_{(p,q)}$) 
	are transversally isotopic.
\end{prop}

\proof
Let $T$ and $T'$ be tori standardly embedded in $M$ on which $K$ and $K'$,  
respectively, sit.

\begin{lem}\label{nonsingular}
	If the self-linking number of $K$ is maximal then $T$ may be
	isotoped relative to $K$ so that the characteristic
	foliation on $T$ is nonsingular.
\end{lem}

This lemma and the next are proved in the following section.

\begin{lem}\label{isotopic}
	Two transversal knots on a torus $T$ with nonsingular characteristic foliation
	that are homologous are transversally isotopic, except possibly when there is a closed leaf
	in the foliation isotopic to the transversal knots.
\end{lem}

Our strategy is to isotop $T$ onto $T'$, keeping $K$ and $K'$ transverse to $\xi$, 
so that $K$ and $K'$ are
homologous, and thus transversally isotopic.  We now show that $T$
can be isotoped into a standard form keeping $K$ transverse (and similarly for
$K'$ and $T'$ without further mention).  Let $V$ be 
the solid torus that $T$ bounds (recall we are choosing $V$ so that $p>q$). 
Let $D_\mu$ and $D_\lambda$ be the disk that $\mu$ and $\lambda$
respective bound. Now observe:

\begin{lem}\label{basis}
	We may take $\mu$ and $\lambda$ to be positive transversal curves and with this
	orientation $\mu,\lambda$ form a positive basis for $T=\partial V$.
\end{lem}

\proof
Clearly we may take $\mu$ and $\lambda$ to be positive transversal knots, for if we could
not then \lm{curveinT} implies that we may isotop one of them to a closed leaf in $T_\xi$ contradicting
the tightness of $\xi$.  Thus we are left to see that 
$\mu, \lambda$ is a positive basis. Assume this is not the case.
By isotoping $T$ slightly we may assume that $T_\xi$ has closed leaf
(indeed if $T_\xi$ does not already have a closed leaf then the isotopy will give an 
intervals worth  of rotation numbers, and hence
some rational rotation numbers, for the return map
induced on $\mu$ by $T_\xi$). Let $C$ be one of these closed leaves and let $n=\lambda\cdot C$
and $m=\mu\cdot C$.  Note $n$ and $m$ are both positive since $\mu$ and $\lambda$ are
positive transversal knots. Since $\mu, \lambda$ is not a positive basis $C$ is an $(n,m)$--torus knot. 
In particular $C$ is a positive torus knot. Moreover, the framing on $C$ induced by $\xi$ is
the same as the framing induced by $T$. Thus $\hbox{tb}(C)=mn$ contradicting \eqn{tb-bound}.
So $\mu,\lambda$ must be a positive basis for $T$.
\endproof

Now let $m=l(\mu)$ and $l=l(\lambda)$ and recall $m,l\leq -1$.

\begin{lem}\label{l-formula}
	If $\gamma$ is a transversal $(p,q)$ knot on $T$ (with nonsingular 
	characteristic foliation)  then 
	\begin{equation}\label{eqn:l-formula}
		l(\gamma)=pm+ql+pq.
	\end{equation}
\end{lem} 

\proof
Let $v$ be a section of $\xi$ over an open 3--ball containing $T$ and its meridional and longitudinal
disks. If $C$ is a curve on $T$ then define
$f(C)$ to be the framing of $\xi$ over $C$ induced by $v$ relative to the framing of $\xi$  
over $C$ induced by $T$.  Note $f$ descends to a map on $H_1(T)$ and $f(A+B)=f(A)+f(B)$
where $A,B\in H_1(T)$. One easily computes $f(\mu)=m$ and $f(\lambda)=l$. Thus 
$f(p\mu+q\lambda)=pm+ql$.  Now for a transversal curve $C$ on $T$ the normal bundle to $C$ can be
identified with $\xi$ thus $f(C)$ differs
from $l(C)$ by the framing induced on $C$ by $T$ relative to the framing induced
on $C$ by its Seifert surface.  So $l(C)=f(C)+pq=pm+ql+pq$. 
\endproof

Thus since $K$ has maximal self-linking number we must have $m=l=-1$. 
Now by \tm{solid_tori} we may find
a contactomorphism from $V$ to $S_f=\{(r,\theta,\phi)\in \R^2\times S^1|
r\leq f(\theta,\phi)\}$ for some positive function $f\co T^2\to \R$, with the 
standard tight contact structure $\hbox{ker}(d\phi+r^2\, d\theta)$. 

Clearly $T=\partial S_f$ may be isotoped to $S_\epsilon=\{(r,\theta,\phi)\in \R^2\times S^1|
r<\epsilon\}$ for arbitrarily small $\epsilon>0$. We now show this isotopy may be done
keeping our knot $K$ transverse to the characteristic foliation. To a foliation on $\partial S_f$ we may
associate a real valued rotation number $r(S_f)$ for the return map on $\mu$ induced
by $(\partial S_f)_\xi$ (see \cite{ml}). 
For a standardly embedded torus this number must be negative since if not then some
nearby torus would have a positive $(r,s)$ torus knot as a closed leaf in its characteristic
foliation violating the Bennequin inequality (as in the proof of Lemma~\ref{basis}). 
So as we isotop $\partial S_f$ to
$\partial S_\epsilon$ we may keep our positive torus knot transverse to the characteristic foliation
by Lemma~\ref{curveinT} (since closed leaves in $(\partial S_f)_\xi$ have slope $r(S_f)$ and
$K$ has positive slope). 
Thus we assume that the solid torus $V$ is contactomorphic to $S_\epsilon$.
If $C$ is the core of $V (=S_\epsilon)$ then it is a transversal unknot with self-linking 
$l(\lambda)=-1$.

Finally, let $V$ and $V'$ be the solid tori associated to the torus knots $K$ and $K'$ and
let $C$ and $C'$ be the cores of $V$ and $V'$.  Now since $C$ and $C'$ are unknots with the
same self-linking number they are transversely isotopic. Thus we may think of $V$ and $V'$ as
neighborhoods of the same transverse curve $C=C'$. From above, $V$ and $V'$ may both be shrunk 
to be arbitrarily small neighborhoods of $C$ keeping $K$ and $K'$ transverse to $\xi$. Hence
we may assume that $V$ and $V'$ both sit in a neighborhood of $C$ which is contactomorphic
to, say, $S_c$ (using the notation from the previous paragraph). By shrinking $V$ and $V'$ 
further we may assume they are the tori $S_\epsilon$ and $S_{\epsilon'}$ inside $S_c$
for some $\epsilon$ and $\epsilon'$. Note that this is not immediately obvious but follows
from the fact that a contactomorphism from the standard model $S_f$ for, say, $V$ to 
$V\subset S_c$ may be constructed to take a neighborhood of the core of $S_f$ to a
neighborhood of the core of $S_c$.
This allows us to finally conclude that we may isotop
$V$ so that $V=V'$. Now since $K$ and $K'$ represent the same homology class on $\partial V$ and
they are both transverse to the foliation we may use Lemma~\ref{isotopic} to transversely 
isotop $K$ to $K'$.
\endproof

\begin{figure}[ht]
	{
\epsfysize=2in\centerline{\relabelbox\small
\epsfbox{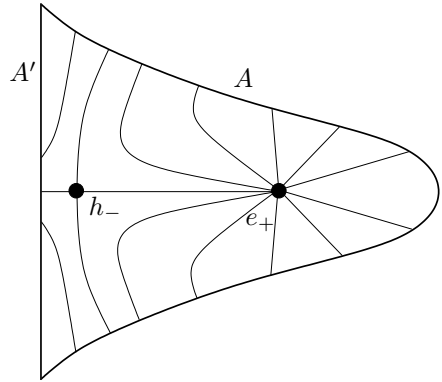}
\relabel {A}{$A$}
\relabel {A'}{$A'$}
\adjustrelabel <-4pt,-1pt> {e}{$e_+$}
\relabel {h}{$h_-$}
\endrelabelbox}}
	\caption{Stabilization disk}
	\label{fig:stab}
\end{figure}

A transversal knot $K$ is called a \dfn{stabilization} of a transversal knot $C$ if
$K=\alpha\cup A$, $C=\alpha\cup A'$ and $A\cup A'$ cobound a disk
with only positive elliptic and negative hyperbolic singularities (eg Figure~\ref{fig:stab}).  
We say $K$ is obtained from $C$ by a \dfn{single stabilization} if $K$ is a stabilization
of $C$ and  $l(K)=l(C)-2$ (ie, the disk that $A\cup A'$ cobound is the one shown in 
Figure~\ref{fig:stab}).
The key observation concerning stabilizations is the following:

\begin{thm}\label{tm-stab}
	If the transversal knots $K$ and $K'$ are single stabilizations of transversal knots
	$C$ and $C'$ then $K$ is transversely isotopic to $K'$ if $C$ is transversely 
	isotopic to $C'$.
\end{thm}

This theorem will be proved in Section~\ref{stab}.  The proof of \tm{main} is completed by an inductive
argument using the following observation.

\begin{lem}\label{canstab}
	If $K$ is a positive transversal $(p,q)$--torus knot and $l(K)<\overline{l}_{(p,q)}$ then $K$
	is a single stabilization of a $(p,q)$--torus knot with larger self-linking number.
\end{lem}

The proof of this lemma will be given in the next section following the proof of \lm{nonsingular}.

% ************************************************************************
\section{Characteristic foliations on tori}\label{tori}
% ************************************************************************

In this section we prove various results stated in Section~\ref{sec:main} related to 
foliations on tori.
Let $T$ be a standardly embedded torus in $M^3$ and $K$ a positive $(p,q)$--torus
knot on $T$ that is transverse to a tight contact structure $\xi$. We are now ready to prove:

\proof[\bf Lemma~\ref{nonsingular}]
{\sl
	If the self-linking number of $K$ is maximal then $T$ may be
	isotoped relative to $K$ so that the characteristic
	foliation on $T$ is nonsingular.
}

\proof
Begin by isotoping $T$ relative to $K$ so that the number of singularities 
in $T_\xi$ is minimal. Any singularities that are left must occur in pairs:
a positive (negative) hyperbolic $h$ and elliptic $e$ point connected by a stable
(unstable) manifold $c$.  Moreover, since $h$ and $e$ cannot be canceled without
moving $K$ we must have $c\cap K\not=\emptyset$.

Now $T\setminus K$ is an annulus $A$ with the characteristic foliation flowing
out of one boundary component and flowing in the other. Let $c'$ be the component
of $c$ connected to $h$ in $A$.  
We can have no periodic orbits in $A$ since such an orbit
would be a Legendrian $(p,q)$--torus knot with Thurston--Bennequin invariant
$pq$ contradicting \eqn{tb-bound}.
Thus the other stable (unstable) manifold $c''$ of $h$ will have to enter (exit) $A$ through the same
boundary component. The manifolds $c'$ and $c''$ separate off a disk $D$ from $A$.
We may use $D\subset T$ to push the arc $K\cap D$ across $D$ to obtain another
transverse $(p,q)$--torus knot $K'$.  It is not hard to show that
$K$ is a stabilization of $K'$. In particular $l(K')>l(K)$, contradicting
the maximality of $l(K)$.  Thus we could have not have had any singularities 
left after our initial isotopy.
\endproof

The above proof provides some insight into \lm{canstab}. Recall:

\proof[{\bf \lm{canstab}}]
{\sl
	If $K$ is a positive transversal $(p,q)$--torus knot with and $l(K)<\overline{l}_{(p,q)}$ then $K$
	is a single stabilization of a $(p,q)$--torus knot with larger self-linking number.
}

\proof
We begin by noting that if $K$ is a stabilization of another transversal knot then it
is also a single stabilization of some transversal knot. Thus we just demonstrate that 
$K$ is a stabilization of some transversal knot.

From the above proof it is clear that if we cannot eliminate all the singularities in
the characteristic foliation of the torus $T$ on which $K$ sits then there is a disk on 
the torus which exhibits $K$ as a stabilization. 

If we can remove all the singularities from $T$ then by \lm{l-formula} we know
that the self-linking number of, say, the meridian $\mu$ is less than $-1$. Thus $\mu$ bounds a 
disk $D_\mu$ containing only positive elliptic and at least one negative hyperbolic singularity.  
To form a positive transversal torus knot $K''$ we can take $p$ copies of the meridian $\mu$ and 
$q$ copies of the longitude $\lambda$
and ``add'' them (ie, resolve all the intersection points keeping the curve 
transverse to the characteristic foliation).  
This will produce a transversal knot on $T$ isotopic to $K$  thus transversely
isotopic.  Moreover, we may use the graph of singularities on $D_\mu$ to show that $K''$, and hence $K$, 
is a stabilization.  
\endproof

We end this section by establishing (a more general version of)
Lemma~\ref{isotopic}.
\begin{lem}
	Suppose that $\mathcal{F}$ is a nonsingular foliation on a torus
	$T$ and $\gamma$ and $\gamma'$ are two simple closed curves on $T.$
	If $\gamma$ and $\gamma'$ are homologous and transverse to $\mathcal{F}$
	then they are isotopic through simple closed curves transverse to  
	$\mathcal{F},$ except possibly if $\mathcal{F}$ has a closed leaf 
	isotopic to $\gamma.$
\end{lem}

\proof
We first note that if $\gamma$ and $\gamma'$ are disjoint and there are not closed
leaves isotopic to them then the annulus that they cobound will provide the desired transverse
isotopy. Thus we are left to show that we can make $\gamma$ and $\gamma'$ disjoint.
We begin by isotoping them so they intersect transversely. Now assume we have transversely 
isotoped them so that the number of their intersection points is minimal.  We wish to show 
this number is zero.  Suppose not, then there are an even number of intersection points
(since homologically their intersection is zero).  

\begin{figure}[ht]
	{\epsfysize=2in\centerline{\epsfbox{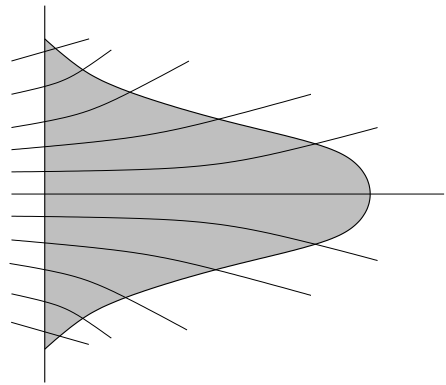}}}
	\caption{Foliation on $D$}
	\label{fig:diskfol}
\end{figure}

Using a standard innermost arc argument we may find a disk $D\subset T$  such that
$\partial D$ consists of two arcs, one a subarc of $\gamma$ the other a subarc of $\gamma'$.
We can use the disk $D$ to guide a transverse isotopy of $\gamma'$ that will decrease the
number of intersections of $\gamma$ and $\gamma'$ contradicting our assumption of minimality.
To see this, note that the local orientability of the foliation implies that we can define
a winding number of $\mathcal{F}$ around $\partial D$. Moreover since $\partial D$ is contractible
and the foliation is nonsingular this winding number must be zero. Thus the foliation on 
$D$ must be diffeomorphic to the one shown in Figure~\ref{fig:diskfol} where the desired isotopy is apparent.
\endproof

% ************************************************************************
\section{Stabilizations of transversal knots}\label{stab}
% ************************************************************************

The main goal of this section is to prove \tm{tm-stab}:

\proof[{\bf \tm{tm-stab}}]
{\sl
	If the transversal knots $K$ and $K'$ are single stabilizations of transversal knots
	$C$ and $C'$ then $K$ is transversely isotopic to $K'$ if $C$ is transversely 
	isotopic to $C'$.
}

\proof
Since $C$ and $C'$ are transversely isotopic we can assume that $C=C'$. 
Let $D$ and $D'$ be the disks that exhibit $K$ and $K'$ as stabilizations of
$C$.  Let $e,h$ and $e',h'$ be the elliptic/hyperbolic pairs on $D$ and $D'$.
Finally, let $\alpha$ and $\alpha'$ be the Legendrian arcs formed by the (closure of the) 
union of stable manifolds of $h$ and $h'$.  Using the characteristic foliation 
on $D$ we may transversely isotop $K\setminus C$ to lie arbitrarily close to $\alpha$
(and similarly for $K'$ and $\alpha'$). We are thus done by the following simple lemmas.

\begin{lem}
	There is a contact isotopy preserving $C$ taking $\alpha\cap C$ to 
	$\alpha'\cap C$.
\end{lem} 

Working in a standard model for a transverse curve
this lemma is quite simple to establish.
Thus we may assume that $\alpha$ and $\alpha'$ both touch $C$ at the same point.
\begin{lem}
	There is a contact isotopy preserving $C$ taking $\alpha$ to $\alpha'$.
\end{lem}
Once again one can use a Darboux chart to check this lemma (for some details see \cite{ef}).
\begin{lem}
	Any two single stabilizations of $C$ along a fixed Legendrian arc are transversely 
	isotopic.
\end{lem}
With this lemma our proof of Theorem~\ref{tm-stab} is complete.
\endproof

We now observe that using \tm{tm-stab} we may reprove Eliashberg's result concerning 
transversal unknots. The reader should note that this ``new proof'' is largely just
a reordering/rewording of Eliashberg's  proof.

\begin{thm}
	Two transversal unknots are transversally isotopic if and only if they
	have the same self-linking number.
\end{thm}

\proof
Using \tm{tm-stab} we only need to prove that two transversal unknots with self-linking
number $-1$ are transversally isotopic, since by looking at the characteristic foliation
on a Seifert disk it is clear that a transversal unknot with self-linking number less
than $-1$ is a single stabilization of another unknot.  But given a transversal unknot with self-linking
number $-1$ we may find a disk that it bounds with precisely one positive elliptic 
singularity in its characteristic foliation. Using the characteristic foliation on the disk 
the unknot may be transversely isotoped into an arbitrarily small neighborhood of the elliptic 
point. Thus given two such knots we may now find a contact isotopy of taking the elliptic point on
one of the Seifert disks to the elliptic point on the other. Since the Seifert disks are tangent
at their respective elliptic points we may arrange that they agree in a neighborhood of the 
elliptic points. Now by shrinking the Seifert disks more we may assume that both unknots sit on the
same disk. It is now a simple matter to transversely isotop one unknot to the other.
\qed

% ************************************************************************
\section{Concluding remarks and questions}
% ************************************************************************

We would like to note that many of the techniques in this paper work for negative torus knots
as well (though the proofs above do not always indicate this). There are two places where we cannot
make the above proofs work for negative torus knots, they are:
\begin{itemize}
\item From Equation~\ref{eqn:l-formula} we cannot conclude that the self-linking numbers of	
	$\mu$ and $\lambda$ are $-1$ when $l(K_{(p,q)})$ is maximal as we could for positive torus
	knots.
\item We cannot always conclude that a negative torus knot with self-linking less than maximal
	is a stabilization.
\end{itemize}
Despite these difficulties we conjecture that negative torus knots are also determined by their
self-linking number.

Let $S=S^1\times D^2$ and let $K$ be a $(p,q)$--curve on
the boundary of $S$. Now if $C$ is a null homologous knot in a three manifold $M$ then 
let $f\co S\to N$ be a diffeomorphism from $S$ to a neighborhood $N$ of $C$ in $M$ taking $S^1\times\{\hbox{point}\}$
to a longitude for $C$.  We now define the \dfn{$(p,q)$--cable of $C$} to be the knot $f(K)$. 
\begin{quest}\label{conj}
	If $\mathcal{C}$ is the class of topological knots whose transversal realizations
	are determined up to transversal isotopy by their self-linking number, then is
	$\mathcal{C}$ closed under cablings?
\end{quest}
Eliashberg's 
Theorem~\ref{unknots} says that the unknot $U$ is in $\mathcal{C}$. Our main Theorem~\ref{main} says
that any positive cable of the unknot is in $\mathcal{C}$. 
This provides the first bit of evidence that the answer to the question might be YES, at least for
``suitably positive'' cablings.

Given a knot type one might hope,
using the observation on stabilizations in this paper, to prove that transversal knots in this
knot type are determined by their self-linking number as follows:  First establishing that there is a 
unique transversal knot in this knot type with maximal self-linking number. Then showing that
any transversal knot in this knot type that does not have maximal self-linking number is a stabilization.
The second part of this program is of independent interest so we ask the
following question:
\begin{quest}
	Are all transversal knots not realizing the maximal self-linking number of their knot type
	stabilizations of other transversal knots?
\end{quest}
It would be somewhat surprising if the answer to this question is YES in complete generality
but understanding when the answer is YES and when and why it is NO should provide insight
into the structure of transversal knots.

We end by mentioning that the techniques in this paper also seem to shed light on Legendrian
torus knots. It seems quite likely that their isotopy class may be determined by their 
Thurston--Bennequin invariant
and rotation number. We hope to return to this question in a future paper.

\end{document}